\documentclass[12pt]{amsart}
\usepackage{amsmath,amsthm,epsfig,amssymb}
\usepackage{amsthm,amsfonts,amssymb,amscd}
\usepackage[all]{xy}
\newtheorem{thm}[subsection]{Theorem}

\newtheorem{lem}[subsection]{Lemma}

\newtheorem{rem}[subsection]{Remark}
\theoremstyle{definition}

\newtheorem{proposition-definition}[subsection]{Proposition-Definition}

\newcommand{\CC}{{\mathbb C}}

\newcommand{\ZZ}{{\mathbb Z}}

\newcommand{\PP}{{\mathbb P}}

\newcommand{\OOO}{{\mathcal O}}

\author{A. El Mazouni}

\address{Laboratoire de Math\'ematiques de Lens EA 2462
Facult\'e des Sciences Jean Perrin
Rue Jean Souvraz, SP18
F-62307 LENS  Cedex France}
\email{mazouni@euler.univ-artois.fr}

\author{D.S. Nagaraj}

\address{Institute of Mathematical Sciences C.I.T. campus, Taramani, 
chennai 600113,India}
\email{dsn@imsc.res.in}

\title{Tangent bundle of $\PP^2$ and morphism from $\PP^2$ to $\text{Gr}(2,
\CC^{4})$}
\date{}
\begin{document}

\begin{abstract}  
In this note we study the image of $\PP^2$ in $\text{Gr}(2, \CC^{4})$ given by tangent bundle of
$\PP^2. $ We show that there is component $\mathcal{H}$ of the Hibert scheme of surfaces in 
$\text{Gr}(2, \CC^{4})$ with no point of it  corresponds to a smooth surface.
\end{abstract}

\maketitle

\section{Introduction} \setcounter{page}{1}
 Let $\PP^2 $ denote the projective plane over the field of complex numbers 
 $\CC$  and $\text{Gr}(2, \CC^{4})$ 
 Grassman variety of two dimensional quotients of the vector space $\CC^4.$
 
The aim of this paper is to study the image of $\PP^2$ by  
non constant morphisms $\PP^2 \to \text{Gr}(2, \CC^{4})$ obtained by tangent 
bundle $T_{\PP^2}$ of $\PP^2.$  The bundle $T_{\PP^2}$ is generated by
sections and hence it is generated by four($=\text{rank}(T_{\PP^2})
+\text{dim}(\PP^2)$) independent
global sections. Any set $S$ of four independent generating sections
of $T_{\PP^2}$ defines a morphism  
$$ \phi_S: \PP^2 \to \text{Gr}(2, \CC^{4}),$$ 
such that the $\phi_S^*(Q) = T_{\PP^2},$ where $Q$ is the universal rank two
quotient bundle on $\text{Gr}(2, \CC^{4}):$
$$\mathcal{O}_{\text{Gr}(2, \CC^{4})}^4 \to Q \to 0.$$

According to a result of Tango \cite{Ta}, if 
$ \phi: \PP^2 \to Gr(2, \CC^{4})$ is an imbedding 
then the pair of Chern
classes $(c_1(\phi^*(Q)), c_2(\phi^*(Q)))$ is 
equal to one of the following
pairs: $(H,0), (H,H^2), (2H,H^2)$ or $(2H,3H^2),$ where $H$ is the ample 
generator of $\text{H}^2(\PP^2, \mathbb{Z}).$ 

Since $c_1(T_{\PP^2}) = 3H$ and 
$c_2(T_{\PP^2}) = 3H^2,$ the morphism $\phi_S$ defined by a
set $S$ of four independent generating sections
of $T_{\PP^2}$ is not an imbedding. Thus, it is natural to ask,
does there exists a  set $S$ of four independent generating sections
of $T_{\PP^2}$ for which $\phi_S$ is generically injective?
In this direction we have the following (Theorem \ref{thm1}):

\begin{thm}
 For  general choice of an ordered set $S$ of four independent generating 
 sections of $T_{\PP^2}$ the morphism 
 $$\phi_S:   \PP^2 \to {\rm Gr}(2, \CC^{4})$$ \begin{footnotesize}              
                               \end{footnotesize}
 is generically injective.
\end{thm}

We also, show that in fact one can find an ordered set of generators $S$
of $T_{\PP^2}$ the morphism is an immersion i.e., the morphism induces
an injection on all the tangent spaces.

As by product of our result we obtain the following (Theorem \ref{hilb}):

\begin{thm}
There is an irreducible component $\mathcal{H}$ of the Hilbert scheme of surfaces in ${\rm Gr}(2, \CC^{4})$
no point  which corresponds  to a smooth surface. 
\end{thm}

\section{The Tangent bundle of $\PP^2.$}

The tangent bundle of $\PP^2$ fits in an exact sequence called the ``Euler
sequence'':
\begin{equation}\label{eq1}
0 \to \mathcal{O}_{\PP^2} \to \mathcal{O}_{\PP^2}(1)^{3} \to 
T_{\PP^n} \to 0. 
\end{equation}
This exact sequence together with the fact $\text{H}^1(\mathcal{O}_{\PP^n})= 0,$
implies that $\text{dim}\text{H}^0(T_{\PP^2}) = 8,$ where $\text{H}^i$
denotes the $i$ th sheaf cohomology group. Since the rank two bundle
$T_{\PP^2}$ on $\PP^2$ is ample and generated by sections, a minimal
generating set of independent sections has cardinality four. 
Any set $S$ of four independent generators of $T_{\PP^2}$ gives to an exact
sequence:
$$ 0\to E_S \to \mathcal{O}_{\PP^2}^4 \to T_{\PP^2} \to 0. $$
This in turn corresponds to a morphism 
$\phi_S :\PP^2 \to \text{Gr}(2,\CC^{4}),$
where $\phi_S(x) = \{ \mathbb{C}^4 = \mathcal{O}_{\PP^2}^4|_x \to 
T_{\PP^2}|_x \to 0\}.$

\section{The main result}

Note that in the ``Euler sequence'' (\ref{eq1}) the injective map
$$0 \to \mathcal{O}_{\PP^2} \to \mathcal{O}_{\PP^2}(1)^{3} $$
is given by the section $v=(X,Y,Z),$ where $X,Y,Z$ is the standard basis 
of $\text{H}^0(\mathcal{O}_{\PP^n}(1)).$ For a section $v_i$ of 
$\mathcal{O}_{\PP^n}(1)^3$ we denote  $w_i$ the image section of $T_{\PP^2}$
under the surjection
$$ \mathcal{O}_{\PP^2}(1)^{3} \to T_{\PP^n} \to 0$$
in (\ref{eq1}). Let 
$\tilde{S} = (v_1,v_2,v_3,v_4)$ be an ordered set of four linearly
independent sections  of $\mathcal{O}_{\PP^2}(1)^{3}$  and 
$S = (w_1,w_2,w_3,w_4)$ be the corresponding ordered set of sections of
$T_{\PP^2}.$
Clearly,the set $S$ is generating set of independent sections of $T_{\PP^2}$ if
and only if $\tilde{S}\cup \{v\}$ is a generating set of independent sections
of $\mathcal{O}_{\PP^2}(1)^{3}.$

\begin{lem}\label{lem1}
Let 
$v_1=(X,0,0),v_2=(0,Y,0),v_3=(Y,Z,X),v_4=(Z, X, Y)$ 
be four sections of $\mathcal{O}_{\PP^2}(1)^3$ and
$\tilde{S}=(v_i|1\leq i \leq 4).$
Then the ordered set $\tilde{S}$ with $ \{v=(X,Y,Z) \}$ is a generating 
set of independent sections
of $\mathcal{O}_{\PP^2}(1)^{3}.$ Hence the corresponding ordered set of sections
$S =(w_1,w_2,w_3,w_4)$  generate $T_{\PP^2},$ where $w_i$ is the image of
$v_i$ under the map given in exact sequence (\ref{eq1}). 
\end{lem}
{\bf Proof:}  Clearly, $v$ generate a subspace of $\mathcal{O}_{\PP^2}(1)^{3}$
dimension one at every point of $\PP^2.$ Hence $w_i,w_j$ is not independent
at a point $p \in \PP^2$ if and only if the section 
$v_{ij}=v_i\wedge v_j \wedge v$
of $\mathcal{O}_{\PP^2}(3)$ vanishes at $p.$ Thus, if the six independent
sections $\{v_i\wedge v_j \wedge v| 1\leq i<j \leq 4\}$ has no common zero
implies $\tilde{S}$ with $\{v=(X,Y,Z) \}$ is a generating set of independent
sections of $\mathcal{O}_{\PP^2}(1)^3.$ Note that $v_{12}=XYZ,v_{13}=X(Z^2-XY),
v_{14}=X(XZ-Y^2), v_{23}=Y(X^2-YZ),v_{24}=Y(YX-Z^2),
v_{34}=3XYZ-(X^3+Y^3+Z^3).$ It is easy to see that the set 
$\{v_{12}, v_{13},v_{14},v_{23},v_{24},v_{34}\}$ of sections of
$\mathcal{O}_{\PP^2}(3)$ has no common
zero in $\PP^2$ and hence the ordered set
$\tilde{S}$ with $\{v=(X,Y,Z) \}$ is a generating set of independent sections
of $\mathcal{O}_{\PP^2}(1)^{3}.$ Hence the corresponding ordered set of sections
$S =(w_1,w_2,w_3,w_4)$  generate $T_{\PP^2}.$ $\hfill{\Box}$

\begin{lem}\label{lem2}
 Let $f: \PP^2 \to \PP^n$ be a non constant morphism and 
 $f^*(\mathcal{O}_{\PP^n}(1))= \mathcal{O}_{\PP^2}(m).$
 Assume that there exists a linear subspace $W$ of codimension
 two such that $W\cap f(\PP^2)$ consists of exactly
 $m^2$ points. Then the morphism $f$ is 
 generically injective.
\end{lem}

{\bf Proof:} Note as the morphism $f$ is non constant
$f^*(\mathcal{O}_{\PP^n}(1))=\mathcal{O}_{\PP^2}(m)$ with $m>0$ and hence is
ample.  This means $f$ is finite map. Set $r= \text{deg}(f),$ the number of
elements $f^{-1}(f(x))$ for a general $x\in \PP^2.$ If $d$ to be the 
degree of $f(\PP^2)$ in $\PP^n$ then it is easy to see that 
$m^2 =d.r.$
On the 
other hand the assumption,  $W\cap f(\PP^2)$ consists of exactly
 $m^2$ points, implies $d\geq m^2.$ Thus we must have $d=m^2$ and
$r=\text{deg}(f)=1.$ Thus $f$ is generically injective.
$\hfill{\Box}$

\begin{lem}\label{lem3}
With the notations of Lemma(\ref{lem1}), the 
surjection of vector bundles on $\PP^2$
$$ \mathcal{O}_{\PP^2}^4 \to T_{\PP^2}$$
given by $S$ defines a generically injective morphism
$$\phi_S : \PP^2 \to {\rm Gr}(2, \CC^{4}).$$
\end{lem}
{\bf Proof:} Let $p: \text{Gr}(2, \CC^{4}) \to \PP^5$ be the Pluker imbedding
given by the determinant of the universal quotient bundle. Then 
$p\circ \phi_S $ is given by 
$$
\begin{array}{l}
(x;y;z) \mapsto \\
(xyz;x(z^2-xy);x(xz-y^2));y(x^2-yz);y(xy-z^2);\\
3xyz-(x^3+y^3+z^3)).
\end{array}
$$
To prove the map $\phi_S$ is generically injective it is enough to prove the map
$p\circ \phi_S $ is so. Set $(Z_0,\ldots, Z_5)$ as the homogeneous
coordinates of $\PP^5$ and $W$ be the codimension two subspace of $\PP^5$
defined by $Z_0= 0 = Z_5.$ Then   
$W\cap p\circ \phi_S(\PP^2)$ is equal to
$$\{(0,-\omega^i, 1,0,0,0);(0,0,0,-\omega^i,
1,0);(0,\omega^i, 1,\omega^i,1,0)|1\leq i\leq 3 \}, 
$$
where $\omega$ is a primitive cube root of unity. Note that 
$$ (p\circ \phi_S)^*(\mathcal{O}_{\PP^5}(1)) = \mathcal{O}_{\PP^2}(3).$$
 Hence,the
required result follows from Lemma(\ref{lem2}). $\hfill{\Box}$

\begin{rem} 
We show (see Lemma \ref{lem5}) that $p\circ \phi_S$
is an immersion. i.e., the induced linear map on the tangent space
at every point of $\PP^2$ is injective and one to one except finitely
many points.
\end{rem}

Next  we recall the follwing [See,
Lemma(3.13)\cite{ALN}]:
\begin{lem}\label{lem4} Let $X$ and $Y$ be two irreducible projective varieties.
Let
$T$ be an irreducible quasi-projective variety and $t_0 \in T$ be a point. Let 
$$F : X \times T \to Y $$
be a morphism. Assume that $F_t :=F|_{X\times t} : X \to Y$ is finite for all
$t\in T$
and $F_{t_0}$ is a birational onto its image.  Then there is an open subvariety 
$U$ of $T$
such that $t_0 \in U$ and for $t \in U$ the morphism $F_t$ is birational onto
its image.
\end{lem}

{\it Proof:} For the sake of completeness we reproduce the proof here.
Consider
the morphism $G=F\times Id_T : X \times T \to Y \times T.$
Then the assumption $F_t$ is finite implies the morphism $G$ is finite and
proper.
Hence ${\mathcal G}= G_{*}(\OOO_{X\times T})$ is coherent  sheaf of
$\OOO_{Y\times T}$ 
modules. Let $Z \subset Y\times T$ be the subvariety on which the sheaf 
$G_{*}(\OOO_{X\times S})$ is supported.  Then clearly the map $p: Z \to T,$
restriction of the natural
projection, is surjective. The section $1 \in \OOO_{X\times T}$ gives an
inclusion of 
$ {\OOO_{Z}}$ in ${\mathcal G}.$ Let ${\mathcal F} = {\mathcal G}/{\OOO_{Z}}.$
Let
$Z_1\subset Y\times T $ be the subvariety on which the sheaf ${\mathcal F} $
supported. 
Let $q : Z_1 \to T$ be the natural projection and let
$U = \{ t \in T| {\rm dim}{q^{-1}(t)} < {\rm dim}(X) \} $ then we see that by
semi continuity
[See, page 95, Exercise (3.22)  \cite{Ha}],
$U$ is an open subset and is non-empty as $t_0 \in U.$   For $t \in U$ the
morphism 
$F_t$ is an isomorphism on $X\times t - G^{-1}(q^{-1}(t).$  Since $G$ is finite 
$G^{-1}(q^{-1}(t)$ is proper closed subset of $X\times t$ and hence the
morphism 
$F_t$  is birational onto its image.
This proves the Lemma.
$\hfill{\Box}$
 
\begin{thm}\label{thm1}
 For a generic choice of an ordered set $S$ of four independent generating 
 sections of $T_{\PP^2}$ the morphism 
 $$\phi_S:   \PP^2 \to {\rm Gr}(2, \CC^{4})$$ 
 is generically injective.
\end{thm}
 
{\it Proof:} It is easy to see that the ordered set of four sections $S$ 
generating  $T_{\PP^2}$
is an irreducible quasi projective variety. In fact it is an open subvariety
of the affine space $V^4,$ where 
$V=\text{H}^0(T_{\PP^2}). $ The theorem  at once
follows from Lemma(\ref{lem4}), if we show the existence of one $S$
for which $\phi_S$ is generically injective. But the 
 existence of one such $S$ follows from Lemma(\ref{lem3}).
 $\hfill{\Box}$
\section{An example}
The result of the previous section can be used give an example of a
component of a Hilbert Scheme of $Gr(2,\CC^4)$ with out any point corresponding
to a smooth surface.

\begin{lem}\label{lem5}
The morphism 
$p\circ \phi_S: \mathbb{P}^2 \to \mathbb{P}^5$ of Lemma(\ref{lem2}) is an immersion i.e., the induced linear map on the tangent space
at every point of $\PP^2$ is injective. Moreover, $p\circ \phi_S$  one to one except  
$$S_1=\{ (1;0;0),(0;1;0), (0;0;1) \} \mapsto (0;0;0;0;0;1)$$
and 
$$ S_2=\{ (1;1;1),(\omega;\omega^2;1),(\omega^2;\omega;1)\} \mapsto (1;0;0;0;0;0), $$
where $\omega$ is a primitive cube root of unity.
\end{lem}

{\it Proof:} Let $X, Y, Z$ be the homogeneous coordinates functions on $\mathbb{P}^2$ 
and $Z_0, Z_1,Z_2,Z_3,Z_4,Z_5 $ be the homogeneous coordinates functions on $\mathbb{P}^5.$
Clearly under the morphism $p\circ \phi_S: \mathbb{P}^2 \to \mathbb{P}^5$ the set 
$S_1$ maps to $(0;0;0;0;0;1)$ and the set $S_2$ maps to $(1;0;0;0;0;0).$ Note that the lines
$X=0, Y=0,$ and $Z= 0$ mapped to nodal cubics $Z_0=Z_1=Z_2=Z_3^3+Z_4^3-Z_3Z_4Z_5 =0,$ 
$Z_0=Z_3=Z_4=Z_1^3+Z_2^3+Z_1Z_2Z_5 =0,$ and $Z_0=Z_1-Z_4=Z_2+Z_3=Z_1^3+Z_2^3-Z_1Z_2Z_5 =0$
respectively. Thus we can conclude that the morphism $p\circ \phi_S$ is an immersion on these
three lines. On the complement of these lines the morphism $p\circ \phi_S$ can be described as  
$$(x,y) \mapsto (1/y-x,x/y-y,x-y/x,y-1/x,3-x^2/y-y^2/x-1/xy) $$
from $\mathbb{C}^2-\{xy=0\} \to \mathbb{C}^4.$ 
If $(x,y)$ and $(x_1,y_1)$  maps to the same point then we get the following equations:
\begin{equation}\label{eq2}
 1/y-x = 1/{y_1}-x_1
\end{equation}
\begin{equation}\label{eq3}
 x/y-y = {x_1}/{y_1}-y_1
\end{equation}\begin{equation}\label{eq4}
 y/x-x= {y_1}/{x_1}-x_1
\end{equation}\begin{equation}\label{eq5}
 1/x-y = 1/x_1-y_1.
\end{equation}
The equations \ref{eq2} and \ref{eq5} gives us 
$$ \frac{xy-1}{y} = \frac{x_1y_1-1}{y_1}; \,\, \frac{xy-1}{x} = \frac{x_1y_1-1}{x_1} .$$
Since $xy\neq 0$ and $x_1y_1\neq 0$ we see that either
$xy-1 \neq 0 $ and  $(x,y) = (x_1,y_1)$ or $xy-1 = 0 =x_1y_1-1 .$
Hence we get $(p\circ \phi_S)$ is one to one out side the set
$ \{(1,1),(\omega,\omega^2),(\omega^2,\omega)\}$  and this is mapped to $(0,0,0,0,0). $ The 
assertion about the immersion of the given morphism 
$$\mathbb{C}^2-\{xy=0\} \to \mathbb{C}^4$$
can be checked by looking at the two by two minors of the below jacobian matrix of the morphism:
$$
\left( \begin{matrix}
  -1 & 1/y & 1+y/{x^2} & 1/{x^2} & -2x/y+y^2/{x^2}+ 1/{x^2y}\\
  -1/{y^2} & -x/{y^2}-1 & -1/x & 1 & {x^2}/{y^2}-2y/x+ 1/{xy^2}
\end{matrix}
\right)
$$
$\hfill{\Box}$

\begin{thm}\label{hilb}
Let $\mathcal{H}$ be the irreducible component of the Hilbert scheme of ${\rm Gr}(2, \CC^{4})$
containing the point corresponding to the image surface  of the morphism
$$\phi_S : \PP^2 \to {\rm Gr}(2, \CC^{4})$$
of \ref{lem3}. Then no point of $\mathcal{H}$ corresponds to a smooth surface. 
\end{thm}

{\bf Proof:} Let $ p:{\rm Gr}(2, \CC^{4}) \to \PP^5$  be the Pluker imbedding.
Since $(p\circ \phi_S)^*(\mathcal{O}_{\PP^5}(1)) = \mathcal{O}_{\PP^2}(3)$ and by Lemma \ref{lem5} the morphism 
$ p\circ \phi_S $ is an imbedding outside finite set of points. Moreover, general hyperplane section
of $(p\circ \phi_S)(\PP^2)$ in $\PP^5$ is smooth curve of genus one. If a point of the irreducible
$\mathcal{H}$ of corresponds to a smooth surface $Y$ then it has to have the same cohomology class as that of
$(p\circ \phi_S)(\PP^2)$
namely $(3,6) \in \rm{H}^4({\rm Gr}(2, \CC^{4}), \ZZ).$ Also, the general hyperplane section of 
$p(Y)$ has to be a smooth curve of genus one. But according to the classification of smooth
surfaces of type $(3,6)$ in ${\rm Gr}(2, \CC^{4})$ (see, \cite[Theorem 4.2]{GM}) implies that
there are no smooth surface of type $(3,6)$ with hyperplane section a smooth curve of 
genus one. This  contradiction proves that no point of $\mathcal{H}$ corresponds to a smooth surface.
$\hfill{\Box}$

\begin{rem} 
The  component $\mathcal{H}$ of the the Hilbert Scheme in Theorem \ref{hilb}
is reduced irreducible of dimension 23. In fact
 computing the normal sheaf associated to the morphism $\phi_S$ of 
Lemma \ref{lem3} and counting the dimension of space of all such 
morphisms we see that $\mathcal{H}$ is a reduced irreducible of 
dimension 23.
\end{rem}

 {\bf Acknowledgment} This work was supported in part by the Labex CEMPI  (ANR-11-LABX-0007-01).
 The second author thanks university D'artois, Lens and the University of Lille. We thank
 Laytimi Fatima for help during the work.

\end{document}